\documentclass{amsart}
\usepackage{amsmath}
\usepackage{amssymb,latexsym,amstext,amsthm}
\usepackage{verbatim} 
\usepackage{eucal} 
\usepackage{color}
\usepackage[draft]{todonotes}
\usepackage{amsfonts}
\usepackage{mdframed}
\usepackage{lipsum}

\usepackage{pgf,tikz}
\usepackage{mathrsfs}
\usetikzlibrary{arrows}
\definecolor{qqqqff}{rgb}{0.,0.,1.}
\definecolor{xdxdff}{rgb}{0.49019607843137253,0.49019607843137253,1.}

\definecolor{qqqqff}{rgb}{0.,0.,1.}

\theoremstyle{plain}
\newtheorem{thm}[equation]{Theorem}
\newtheorem{pro}[equation]{Proposition}
\newtheorem{cor}[equation]{Corollary}
\newtheorem{lem}[equation]{Lemma}
\theoremstyle{remark}

\newtheorem{ex}[equation]{\bf Example}

\newtheorem{DEF}[equation]{\bf Definition}

\numberwithin{equation}{section}

\begin{document}
\setcounter{page}{1}
\title[Various notions of point derivations in Banach algebras]{Various notions of point derivations in Banach algebras}


\author{S.Chameh and A. Rejali*}
\address
{S.Chameh\\}
\address
{Department of  Pure Mathematics\\ University of Isfahan\\Isfahan, Iran,\\}
 \email{S.chameh@sci.ui.ac.ir}\email{Sakineh.chameh@gmail.com}
\address
{A.Rejali}
\address
{Department of Pure Mathematics\\ University of Isfahan\\Isfahan, Iran,\\}
\email{rejali@sci.ui.ac.ir}
\email{ali.rejali2020@gmail.com }

\thanks{
* Corresponding author.
}
\maketitle
\begin{abstract}
 In this paper, we study several type of point derivations for Banach algebras. We investigate how our definition of point derivations are related to each others. 
\end{abstract}\vspace{2mm}

\section{Introduction and Preliminaries}
\indent One of the most interesting concepts in harmonic analysis is Lipschitz Banach algebra which is studied first by Myers \cite{My}. Unfortunately he couldn't publish his achievements. His interest in Lipschitz algebras was connected with the process of differentiation and it was followed by Sherbert \cite{Sh,She}. He supplied proofs of many of the unproved statements of Myers and in some cases extended his result to more general settings.\\
\indent Johnson \cite{Joh} introduced  and studied some basic properties of vector-valued Lipschitz algebras. There are valuable papers related to  the notion of  vector-valued Lipschitz algebras; for example, see \cite{ABR,Biya,Ra,Esm}. The various notions of amenability of vector-valued Lipschitz algebras  have been introduced and studied. We refer the reader to \cite{AAR,BA,Hu,Ran}, for more information.\\
\indent Derivation plays an essential role in some important branches of physics such as dynamic systems \cite{Mir,Abb}. Another application of derivations in mathematics is the concept of amenability for a Banach algebra which was introduced by Johnson \cite{Jo}. Mirzavaziri and Moslehian introduced and studied $\sigma$-derivation. For more details, see \cite{Mir,Mirz}.\\
\indent A point derivation which is generalization of differentiation at a point, is a linear functional on algebra of functions which satisfies the product rule. They were first studied  on  an algebra by Singer and Warmer \cite{Si}. Sherbert \cite{Sh} looked at the purely algebraic aspects of point derivations. He presented three characterizations of the structure of point derivations on Lipschitz algebra. The first in terms of cluster points of certain sequences of functionals, second in terms of Banach limits and the third in terms of the Stone-Cech compactification of certain space. The first two treatments hold in any metric space $(X,d)$ but for the third it requires $(X,d)$ to be a certain compact space. Myers \cite{My} stated the first two results.\\
\indent Let $(X,d)$ be a metric space and $A$ be a complex Banach algebra. An $A$-valued function $f$ on $X$ is called Lipschitz function, if there exists a constant $K>0$ such that 
$$\| f(x)-f(y)\|_{A}\leq K d(x,y)\hspace*{8mm}(x,y\in X).$$
For $A$-valued function $f$ on $X$, let
$$P_{d,A}(f):=\sup\{\dfrac{\| f(x)-f(y)\|_{A}}{d(x,y)}: x,y\in X \hspace*{2mm} and\hspace*{2mm} x \neq y\},$$
and $$\| f\|_{\infty, A}=\sup\lbrace \| f(x)\|_{A} : x\in X\rbrace.$$
Then $Lip_{d}(X,A)$ is the collection of all bounded $A$-valued functions $f$, which $P_{d,A}(f)<\infty$, which is a Banach algebra under the pointwise multiplication and the norm
 $$\Vert f\Vert_{d,A}:=\Vert f\Vert_{\infty, A}+P_{d,A}(f)\hspace*{8mm}(f\in Lip_{d}(X,A)).$$ \\
\indent It is important to know that  the algebraic properties of these algebras depend on the algebraic properties of $A$ and the topological properties of metric space $(X,d)$, see \cite{ABR}. For example, it is easy to see that $Lip_{d}(X,A)$ is a commutative (unital) Banach algebra if and only if $A$ is a commutative (unital) Banach algebra. Note that, whenever $A=\mathbb{C}$, we denote the classical complex-valued Lipschitz algebra by $Lip_{d}X$, which was first studied by Sherbert, see \cite{She,Sh}.\\
\indent Here, we provide some preliminaries which are required throughout the paper. Let $(X,d)$ be a metric space with at least two elements. A linear functional $D$ on $Lip_{d}X$ is called a point derivation at $x$ if 
$$D(fg)=f(x) Dg+ g(x) Df \hspace*{6mm}( f,g \in Lip_{d}X),$$
 Let $\mathfrak{D}_{x}$ denote the set of all point derivations on $Lip_{d}X$ at $x\in X$. Then $\mathfrak{D}_{x}$ is a linear space. For more details readers may refer to \cite{Sh}. Let $A$ be a Banach algebra and $X$ be a Banach $A$-bimodule. A bounded linear map $d:A\longrightarrow X$ is called a derivation if
$$d(ab)=a.db+(da).b\hspace*{6mm}(a,b\in A),$$
A derivation $ad_{x}:A\longrightarrow X$ is called an inner-derivation, if there exists $x\in X$ such that 
$$ad_{x}(a)=a.x-x.a\hspace*{6mm}(a\in A).$$
The set of all derivations is denoted by $Z^{1}(A,X)$ and the set of all inner-derivation is denoted by $B^{1}(A,X)$. Clearly, $B^{1}(A,X)$ is a (not necessarily closed) subspace of $Z^{1}(A,X)$, see \cite [p.38]{Run} for more details.The character on a commutative Banach algebra $A$ is  the space, consisting of all nonzero multiplicative linear functionals on $A$, denoted by $\Delta(A)$. We equipe $\Delta(A)$ with the Gelfand topology, which is a relative topology on $\Delta(A)\cup\{0\}$ induced by the weak$^\ast$-topology of $A^{\ast}$, where $A^{\ast}$ is the dual space of $A$. If $\Delta(A)$ is not empty, then for every $a\in A$, the mapping $\widehat{a}:\Delta(A)\longrightarrow\mathbb{C}$ is defined by $\widehat{a}(\varphi)=\varphi(a)$ for all $\varphi\in \Delta(A)$. It is called the Gelfand transform of $a$.\\
\indent  We study the notion of point derivation on vector-valued Lipschitz algebras. It should be noted that our results are valid for many Banach algebras of functions, however we are interested in studying Lipschitz Banach algebras. Moreover, we investigate how our definition of point derivation and other definitions are related to each other. In some cases we find out relation, for example Runde's definition \cite[p.38]{Run}, Sherber's definition \cite{Sh}  and Mirzavaziri's definition \cite{Mir}.
\section{some algabraic properties of point derivation}
\indent In this section, we  first briefly discuss some of the basic properties of the point derivation on Banach algebras. Afterward we state our main results which compare different types of derivations with our definition.\\
\indent Throughout this section, we will make the following assumptions:\\
Let $A$ and $B$ be Banach algebras, $Hom(A,B)$ stands for continuous homomorphism from $A$ into $B$. 
\begin{DEF}
Let $\varphi\in \Delta(B)$ and $\psi\in Hom(A,B)$. A bounded linear functional $D$ on a Banach algebra $A$ is called $(\varphi,\psi)$-point derivation on  Banach algebra $A$ if 
$$ D(ab)=\varphi\circ\psi(a) Db+\varphi\circ\psi(b) Da\hspace*{8mm}(a,b\in A).$$
Let us denote the set of all $(\varphi,\psi)$-point derivations on $A$ by $\mathfrak{D}_{(\varphi,\psi)}$.
\end{DEF}
The following example illustrates well our definition of point derivation. 
\begin{ex}
Let $A:=C^{1}[0,1]$ and $B:=C[0,1]$ and $D:A\longrightarrow\mathbb{C}$ defined by $Df=f'(x_{\circ})$, for all $f \in A$ and $x_{\circ}\in[0,1]$ . Then the linear functional $D$ is a $(\varphi,\psi)$-point derivation on $A$, where $\psi=id_{A}$ and $\varphi=\delta_{x_{\circ}}$, $x_{\circ}\in[0,1]$
 \end{ex}
\begin{pro}
 For every $\varphi\in \Delta(B)$ and $\psi\in Hom(A,B)$, $\mathfrak{D}_{(\varphi,\psi)}$ is a $w^{\ast}$- closed subspace of $A^{\ast}$. 
\end{pro}
 \begin{proof}
It is easy to check that $\mathfrak{D}_{(\varphi,\psi)}$ is a linear subspace of $A^{\ast}$. We just prove that $\mathfrak{D}_{(\varphi,\psi)}$ is weak*-closed. Let $D_{\circ}$ be a cluster point of $\mathfrak{D}_{(\varphi,\psi)}$ in weak*-topology. Then there exists a net $(D_{\alpha})_{\alpha}$ in $\mathfrak{D}_{(\varphi,\psi)}$ weak*-convergence to $D_{\circ}$ in weak*-topology.\\
For every $a$ and $b$ in $A$,
 \begin{align}
D_{\circ}(ab)=\lim_{\alpha}D_{\alpha}(ab)=&\lim_{\alpha}(\varphi\circ\psi(a)D_{\alpha}(b)+\varphi\circ\psi(b)D_{\alpha}\nonumber(a))\\\nonumber
=&\lim_{\alpha}\varphi\circ\psi(a)D_{\alpha}(b)+\lim_{\alpha}\varphi\circ\psi(b)D_{\alpha}(a)\\\nonumber
=&\varphi\circ\psi(a)D_{\circ}(b)+\varphi\circ\psi(b)D_{\circ}(a).\nonumber
\end{align}
\end{proof}
\indent Now, we state the following useful Lemma.
\begin{lem}\label{1}
Suppose that $A$ and $B$ are Banach algebras of complex-valued functions defined on a set $X$ such that contain constant functions on $X$, $D\in \mathfrak{D}_{(\varphi,\psi)}$, where $\varphi\in \Delta(B)$ and $\psi\in Hom(A,B)$. Then the following statements hold:\\
(i) For every $f\in A$ and every $n\in\mathbb{N}$, $Df^{n}=n(\varphi\circ\psi(f))^{n-1}Df$;\\
(ii) If $f\in A$ satisfies $f^{2}=f$, then $Df=0$;\\
(iii) If $f\in A$ is constant function, then $Df=0$, for all $D\in \mathfrak{D}_{(\varphi,\psi)}$;\\
(iv) Let $\mathfrak{P}$ denote a family of idempotent elements of $A$ and $\overline{\langle\mathfrak{P}\rangle}^{\Vert.\Vert_{A}}=A$. Then $ \mathfrak{D}_{(\varphi,\psi)}=\{0\}$.
\end{lem}
\begin{proof}
$(i)$ It is easily proved by induction.\\
$(ii)$ The proof falls naturally into two parts. First, consider $\psi(f)= 0$. By using part $(i)$, we have $$Df=Df^{2}=2\varphi\circ\psi(f)Df=0.$$Now, we can suppose that $\psi(f)\neq 0$. If $\varphi\circ\psi(f)= 0$, by a similar argument, one can show that $Df=0$. Let $\varphi\circ\psi(f)\neq 0$ and $\psi(f)\neq 0$. Since $f$ is idempotent, then by using part $(i)$ for $f^{2}$ and $f^{3}$, It is easy to see that $\varphi\circ\psi(f)Df= 0$. It completes the proof.\\$(iii)$ Define $1(x)=1_{\mathbb{C}}$, for all $x\in X$. By using part $(ii)$, $D1=0$. Now, let $f$ be a constant function. Then, $f=\lambda . 1$, for some $\lambda\in{\mathbb{C}}$. Since $D$ is linear, we have $Df=0$.\\ $(iv)$ Let $f\in A=\overline{\langle\mathfrak{P}\rangle}^{\Vert.\Vert_{A}}$. Then there exists a sequence $\{f_{n}\}$ in $\langle\mathfrak{P}\rangle$ converges to $f$. Take $D\in \mathfrak{D}_{(\varphi,\psi)}$ and $g\in \langle{P}\rangle$, we have$$ g=\sum^{n}_{i=1}\lambda_{i}g_{i}\hspace*{6mm} ( \lambda_{i}\in\mathbb{C},g_{i}\in \mathfrak{P},1\leq i\leq n),$$where, $g_{i}^{2}=g_{i}$, $1\leq i\leq n$. By using part $(ii)$, $Dg=0$. Since $D$ is continuous and $Df_{n}=0$, for all $n\in\mathbb{N}$. It follows easily that $Df=0$. 
\end{proof}
\indent It should be noted that the previous proposition expresses if $Df=0$, then it is not necessary that $f$ to be a constant function.
\begin{pro}\label{2}
 Let $A$ and $B$ be Banach algebras of complex-valued functions defined on a set $X$ such that contain constant functions on $X$, $\psi\in Hom(A,B)$ and $\varphi\in \Delta(B)$. If $f\in A$ satisfies a polynomial identity $P(f)=0$, where $P$ is a polynomial in one variable, then $Df=0$ for all $(\varphi,\psi)$-point derivations $D$on $A$.
\end{pro}
\begin{proof}
Since $P(f)=0$, the range of $f$ is a finite set. Let the set $\{\lambda_{1},\lambda_{2},...,\lambda_{n}\}$ be the range of $f$. Put $E_{i}:=\{ x\in X: f(x)=\lambda_{i}\}$, $1\leq i\leq n$. Thus,$f=\sum^{n}_{i=1}\lambda_{i}\chi_{E_{i}}$ where the family $\{E_{i}\}^{n}_{i=1}$ is a partition of $X$.\\For every $j$, $1\leq j\leq n$, define $$Q_{j}(\lambda)=\Pi^{n}_{i=1,i\neq j}\frac{(\lambda-\lambda_{i})}{(\lambda_{j}-\lambda_{i})}\hspace*{6mm}(\lambda\in\mathbb{C}).$$ We have $Q_{j}(\lambda_{i})=\delta_{ij}$, for each $1\leq i,j\leq n$. Since $A$ is an algebra, we have $Q_{j}\circ f\in A$. It is easy to check that $Q_{j}\circ f=\chi_{E_{j}}$, for every $j$, $1\leq j\leq n$. Therefore by using Lemma \ref{1} part $(ii)$, $D\chi_{E_{j}}=0$. This completes the proof.
\end{proof}
 \begin{cor} 
Let $A$ and $B$ be Banach algebras of complex-valued functions defined on a set $X$ such that contain constant functions on $X$, $\psi\in Hom(A,B)$ and $\varphi\in \Delta(B)$. If $F$ denote the linear span of the set $\{ f:f\in A, P(f)=0\}$, where $P$ is a polynomial in one variable such that $F$ is dense in $(A, \Vert.\Vert_{A})$, then $ \mathfrak{D}_{(\varphi,\psi)}=0$. 
\end{cor}
\begin{proof}
Suppose that $f\in A$ and $D\in\mathfrak{D}_{(\varphi,\psi)}$, then there exists a sequence $\{f_{n}\}$ in $F$ such that $\{f_{n}\}$ converges to $f$. For every $n\in\mathbb{N}$, $f_{n}=\sum^{l}_{i=1}\alpha_{k}g_{n,k}$, where $g_{n,k}\in A$ and $P(g_{n,k})=0$, $1\leq k\leq l$. By using Proposition \ref{2}, $Df_{n}=0$, for all $n\in\mathbb{N}$. Since $D$ is continuous, it implies $Df=0$.
\end{proof}
\indent The next proposition asserts for every nonzero point derivation, there exists unique $\psi\in Hom(A,B)$ such that $D\in \mathfrak{D}_{(\varphi,\psi)}$, whenever $B$ is a commutative semisimple Banach algebra.
\begin{pro}
Suppose that  $A$ is a Banach algebra and B is a commutative semisimple Banach algebra. Let $D$ be a non-zero point derivation on A such that $ D\in \mathfrak{D}_{(\varphi,\psi_{1})}\cap \mathfrak{D}_{(\varphi,\psi_{2})}$. Then $\psi_{1}=\psi_{2}$.
\end{pro}
\begin{proof}
Since $D$ is not zero, there exists, $f\in A$ such that $Df\neq0$. We may assume that $\varphi(\psi_{1}(f))=0$ for all $\varphi\in\Delta(B)$. Otherwise  replace $f$ by $f^{2}-\varphi(\psi_{1}(f))f$. Since $B$ is semisimple, we have $\psi_{1}(f)=0$. By using Lemma \ref{1} part $(i)$ for $f^{2}$ and $f^{3}$ and based on our hypothesis $ D\in \mathfrak{D}_{(\varphi,\psi_{1})}\cap \mathfrak{D}_{(\varphi,\psi_{2})}$, we obtain
$$\varphi(\psi_{1}(f)) Df=\varphi(\psi_{2}(f)) Df \hspace{6mm}(\varphi\in \Delta(B)).$$
Since $B$ is semisimple, $\psi_{2}(f)=0$. Let $g\in A$. By easy computation, one can show that 
$$Dfg=\varphi(\psi_{1}(g)) Df=\varphi(\psi_{2}(g)) Df \hspace{6mm}(\varphi\in \Delta(B)).$$  
Hence
 $$\psi_{1}(g)-\psi_{2}(g)\in ker\varphi\hspace*{3mm}(\varphi\in\Delta(B)).$$
 Since $B$ is semisimple, it implies that $\psi_{1}=\psi_{2}$.
\end{proof}
\indent In the same manner we can state the next proposition which is proved in a similar way.
\begin{pro}Let $A$ be a Banach algebra, $\psi\in Hom(A)$ and $D$ be a non-zero point derivation on A such that $ D\in \mathfrak{D}_{(\varphi_{1},\psi)}\cap \mathfrak{D}_{(\varphi_{2},\psi)}$, where $\varphi_{1},\varphi_{2}\in\Delta(A)$. Then $\varphi_{1}=\varphi_{2}$.
\end{pro}
\begin{proof}
Since $D$ is not zero, there exists $f\in A$ such that $Df\neq0$. We can assume that $\varphi_{1}\circ\psi(f)=0$. Otherwise replace $f$ by $f^{2}-\varphi_{1}\circ\psi(f)f$. Since $D\in\mathfrak{D}_{(\varphi_{1},\psi)}\bigcap\mathfrak{D}_{(\varphi_{2},\psi)}$, by using Lemma \ref{1} part $(ii)$ for $f^{2}$, we obtain $\varphi_{2}\circ\psi(f)=0$. Let $g\in A$. Then $$Dfg=\varphi(\psi_{1}(g)) Df=\varphi(\psi_{2}(g)) Df \hspace{6mm}(\psi\in Hom(A)).$$Now, by simple argument we conclude $$(\varphi_{1}-\varphi_{2})(\psi(g))=0\hspace*{6mm}(\psi\in Hom(A),g\in A).$$Put $\psi=id_{A}$ . The proof is complete.\end{proof}
\indent In what follows, we state that some theorems yield information about the relation between our definition and Runde's definition  \cite[p.38]{Run}.
\begin{thm}\label{3}
Suppose that $A$ is a Banach subalgebra of a semisimple Banach algebra $B$ and $D:A\longrightarrow B$ is a bounded linear operator. Then the following statements are equivalent:\\
(i) $D$ is a derivation;\\
(ii) $\varphi\circ D$ is $(\varphi,id_{A})$-point derivation, for all $\varphi\in\Delta (B)$. 
\end{thm}
\begin{proof}
 Sufficiently clear. For necessity, let $a,b\in A$ and $\varphi\in\Delta (B)$. Then
\begin{align}\nonumber
\varphi\circ D(ab)=&\varphi\circ id_{A}(a) \varphi\circ Db+\varphi\circ id_{A}(b) \varphi\circ Da\\\nonumber
=&\varphi(a)\varphi\circ Db+\varphi(b) \varphi\circ Da\\\nonumber
=&\varphi(a.Db+b.Da)\nonumber
\end{align}
Since $B$ is semisimple, It follows that $D(ab)=a.Db+b.Da$, for all $a,b\in A$.
\end{proof}
\begin{DEF}
Let $A$ be Banach algebra, $X$ be $A$-bimodule Banach algebra. The family of  all non-zero complex multiplicative linear functional on $X$ which is $A$-bimodule homomorphism  is denoted by $\Delta_{A}(X)$.\\
We call $X$ is semisimple  $A$-bimodule, if $R( A,X)=\bigcap_{\varphi\in\Delta_{A}(X) }ker\varphi=\{0\}$. In particular, $A=\mathbb{C}$ then $\Delta_{A}(X)=\Delta(X)$ and $R(A,X)=R(X)=\bigcap_{\varphi\in\Delta(X) }ker\varphi$.
\end{DEF}
\begin{thm}
Suppose that $A$ is Banach algebra, $X$ is semisimple and commutative $A$-bimodule Banach algebra. If $D:A\longrightarrow X$ is bounded linear operator, then the following statements equivalent:\\
(i) $D$ is derivation;\\
(ii) $\varphi\circ D$ is the derivation for all $\varphi\in\Delta_{A}(X)$.
\end{thm}
\begin{proof} 
Sufficiently, is straightforward. Conversely, for every $a,b\in A$
$$\varphi\circ D(ab)=a.\varphi\circ Db+\varphi\circ Da .b=\varphi(a.Db+Da.b).$$
Since $X$ is semisimple $A$-bimodule, we obtain 
$$D(ab)=a.Db+Da.b \hspace*{6mm}( a,b\in A).$$
\end{proof}
Now, we are ready to illustrate another fact about our definition and Runde's definition \cite[p.38]{Run}.
\begin{thm}
 Suppose that $A$ is a Banach algebra, $X$ is a $A$-bimodule Banach algebra and $D:A\longrightarrow X$ is bounded linear operator. If there is $\varphi\in\Delta_{A}(X)$ which is bijective, then the following statements equivalent:\\
(i)  $D$ is inner-derivation.\\
 (ii)  $\varphi\circ D$ is inner-derivation.
\end{thm}
\begin{proof}
$(i) \Rightarrow (ii)$ Clearly, $\varphi\circ D$ is a derivation. Since $D$ is inner-derivation, there exists $x_{\circ }\in X$ such that $D=ad_{x_{\circ }}$. Put $a_{\circ }=\varphi(x_{\circ })$. By easy computation, one can show that $\varphi\circ D=ad_{a_{\circ }}$.\\
$(ii)\Rightarrow (i)$ Since $\varphi\circ D$ is inner, there is $a\in A$ such that $\varphi\circ D=ad_{a}$. On the other hand, $\varphi$ is epimorphism, therefore there exists $x\in X$ such that $\varphi(x)=a$. Now for all $b\in A$, we have 
$$\varphi\circ D(b)= ad_{a}(b)=b.a-a.b= b.\varphi(x)-\varphi(x).b= \varphi(b.x-x.b)$$
Since $ \varphi $ is injective, we have 
$$Db=b.x-x.b\hspace*{6mm}(b\in A).$$
\end{proof}

\begin{pro}
Let $A$ be a Banach algebra, $X$ be a $A$-bimodule Banach algebra and $D:A\longrightarrow X$ be an inner-derivation. Then the following statements hold:\\
(i) $\varphi\circ D= 0$ for each $\varphi\in\Delta_{A}(X)$.\\
(ii) If $X$ is semisimple $A$-bimodule, then $D= 0$.
\end{pro}
\begin{proof}
$(i)$ It is standard.\\
$(ii)$ According part $(i)$, for all $a\in A$ and $\varphi\in\Delta_{A}(X)$, we have $Da\in ker\varphi$. Since $X$ is semisimple, it follows that $Da=0$, for all $a\in A$ .
\end{proof}
\indent Let us mention another result that may be proved in much the same way as Theorem \ref{3}. For convenience, we repeat the relevant material from \cite{Mir}.\\
Suppose that $A$ is a Banach algebra, $\psi\in Hom(A)$ and $X$ is a Banach $A$-bimodule. Then \\
$(i)$ A linear operator $D:A\longrightarrow X$ is $\psi$-derivation if satisfies
 $$D(ab)=D(a)\psi(b)+\psi(a)D(b)\hspace*{6mm}(a,b\in A).$$
$(ii)$ A $\psi$-derivation $D$ is $\psi$-inner derivation if there is a $x\in X$ such that
$$D(a)=\psi(a).x-x.\psi(a)\hspace*{6mm}(a\in A).$$
\indent In fact, the next theorem yields information about relation between our definition and Mirzavaziri's definition \cite{Mir}.
\begin{thm}
Let $A$ be Banach algebra, $X$  be commutative semisimple $A$-bimodule Banach algebra and $D:A\longrightarrow X$ is bounded linear operator. Then the following statement equivalent:\\
(i) $D$ is $\psi$-derivation;\\
(ii) $\varphi\circ D$ is $\varphi\circ\psi$-derivation, forall $\varphi\in\Delta (X)$;\\
(iii) $ \varphi\circ D$ is $(\varphi,\psi)$-point derivation, for all $\varphi\in\Delta (X)$.
\end{thm}
\begin{proof}
$(i)\Rightarrow (ii)$ By using definition of $\varphi\circ\psi$-derivation, it is straightforward.\\
$(ii)\Rightarrow (iii)$ It is clear.\\
$(iii)\Rightarrow (i)$ Let $a,b\in A$. Then
\begin{align}\nonumber
\varphi\circ D(ab)=&\varphi\circ\psi(a) \varphi\circ D(b)+\varphi\circ\psi(b) \varphi\circ D(a)\\\nonumber
=& \varphi(\psi(a)Db)+\varphi(\psi(b)Da)\\\nonumber
=&\varphi(\psi(a)Db+\psi(b)Da).\nonumber
\end{align}
Since $X$ is commutative semisimple $A$-bimodule Banach algebra, we have 
$$D(ab)=\psi(a)Db+\psi(b)Da\hspace*{6mm}(a,b\in A).$$
\end{proof}
\indent The remainder of this section will be devoted to generalizing the concept of point derivation at $x$, which was introduced by Sherbert \cite{Sh}.
\begin{DEF}
Let $(X,d)$ be a metric space, $E$ be a Banach algebra and $A$ be an algebra of $E$-valued functions. We call the linear functional $D$ on $A$ , is a $(x,\varphi)$-point derivation if 
$$D(fg)= \varphi(f(x))Dg+\varphi(g(x))Df\hspace*{6mm}(f,g\in A).$$
where  $x\in X$ and $\varphi\in\Delta(E)$.
  The set of all $(x,\varphi)$-point derivations on $A$ is denoted by $\mathfrak{D}_{(x,\varphi)}$.
\end{DEF}
\begin{lem} 
 Suppose that $(X,d)$ is a metric space, $E$ is a Banach algebra and $A=Lip_{d}(X,E)$. Also suppose that  $P$ is polynomial in one variable with $P(f)=0$ for some $f\in A$. If $D:Lip_{d}X\longrightarrow\mathbb{C}$ is point derivation based on Sherbert's definition, then $D(\varphi\circ f)=0$, for all $\varphi\in\Delta(E)$.
\end{lem}
 \begin{proof}
Let $P(z)=\sum_{k=1}^{n} c_{k}z^{k}$. Since $P(f)=0$ we have $(P(f))(x)=0$ for each $x\in X$. For $\varphi\in\Delta(E)$, we have
  $$\sum_{k=1}^{n} c_{k}(\varphi(f(x)))^{k}=0\hspace*{3mm}(x\in X).$$
  Note that we have actually shown that the range of $\varphi\circ f$ is finite. Thus there exist $x_{1},x_{2},...,x_{n}\in X$ such that $\sum_{k=1}^{n} c_{k}(\varphi(f(x_{j})))^{k}=0$,$1\leq j\leq n$. Define $\lambda_{j}=\varphi\circ f(x_{j})$ and $E_{j}=\{ x\in X:\varphi\circ f(x)=\lambda_{j}\}$, for each $j$, $1\leq j\leq n$. It follows that $\varphi\circ f=\sum_{j=1}^{n}\lambda_{j}\chi_{E_{j}}$. By using \cite [Lemma 2.1]{BA}, $\varphi\circ f\in Lip_{d}X$. Now, according to Proposition \ref{2}, we conclude $D(\varphi\circ f)=0$ for every  $\varphi\in\Delta(E)$.
 \end{proof}
\indent  In the next theorem which is our main result of this section, we show that  how our definition of point derivation and Sherbert's \cite{Sh} are related to each other.
\begin{thm}
 Let $(X,d)$ be a metric space, $A$ be a Banach algebra and $\varphi\in\Delta(A)$. If $D:Lip_{d}X\longrightarrow\mathbb{C}$ is a bounded linear functional, then the following statement are equivalent:\\
 (i) $D$ is a point derivation at $x$;\\
 (ii) $D_{\varphi}: Lip_{d}(X,A)\longrightarrow\mathbb{C}$ defined by $D_{\varphi}(f)=D(\varphi\circ f)$ is $(x,\varphi)$-point derivation.
\end{thm} 
\begin{proof}
$(i)\Rightarrow (ii)$.Let  $f,g\in Lip_{d}(X,A)$, $\varphi\in\Delta(A)$. Then
\begin{align}\nonumber
D_{\varphi}(fg)=D(\varphi\circ fg)=&D((\varphi\circ f)(\varphi\circ g))\\\nonumber
=&\varphi\circ f(x) D(\varphi\circ g)+ \varphi\circ g(x)D(\varphi\circ f)\\\nonumber
=& \varphi(f(x)) D_{\varphi}g+ \varphi(g(x)) D_{\varphi} f.\nonumber
\end{align}
$(ii)\Rightarrow (i)$. Choose $\varphi\in\Delta(A)$. Since $\varphi$ is nonzero there is $z\in A$ such that $\varphi(z)\neq 0$. Let $f_{\circ},g_{\circ}\in Lip_{d}X$. Define $f(x)=\frac{z}{\varphi(z)}f_{\circ}(x)$ and $g(x)=\frac{z}{\varphi(z)}g_{\circ}(x)$ for all $x\in X$. Therefore, $f,g\in Lip_{d}(X,A)$ such that $g_{\circ}=\varphi\circ g$ and $f_{\circ}=\varphi\circ f$. Since $D_{\varphi}$ is $(x,\varphi)$-point derivation, we conclude
\begin{align}\nonumber
D( f_{\circ}g_{\circ})=&D((\varphi\circ f)(\varphi\circ g))\\\nonumber
=&D(\varphi\circ fg)\\\nonumber
=&D_{\varphi}(fg)\\\nonumber
=& \varphi\circ f(x)D_{\varphi} g+ \varphi\circ g(x)D_{\varphi} f\\\nonumber
=&\varphi\circ f(x)D(\varphi\circ g)+ \varphi\circ g(x) D(\varphi\circ f)\\\nonumber
=&f_{\circ}(x)Dg_{\circ}+g_{\circ}(x)Df_{\circ}\nonumber
\end{align}
It completes the proof.
\end{proof}

\end{document}